\documentclass[12pt]{amsart}
\usepackage{verbatim}
\usepackage{amsmath}
\usepackage{amssymb}
\usepackage{amscd}

\numberwithin{equation}{section}

\newtheorem{theorem}[equation]{Theorem}
\newtheorem{lemma}[equation]{Lemma}

\newtheorem{proposition}[equation]{Proposition}
\newtheorem{corollary}[equation]{Corollary}

\newtheorem*{theorem*}{Theorem}

\theoremstyle{definition}

\theoremstyle{remark}
\newtheorem{remark}[equation]{Remark}

\newcommand{\la}{{\langle}}
\newcommand{\ra}{{\rangle}}

\newcommand{\X}{{\mathcal X}}

\newcommand{\C}{{\mathcal C}}

\newcommand{\M}{{\mathcal M}}
\newcommand{\N}{{\mathcal N}}

\def\CC{{\mathbb C}}

\def\NN{{\mathbb N}}
\def\RR{{\mathbb R}}
\def\ZZ{{\mathbb Z}}

\newcommand{\Div}{\operatorname{Div}}
\newcommand{\divisor}{\operatorname{div}}
\newcommand{\Jac}{\operatorname{Jac}}
\newcommand{\Prin}{\operatorname{Prin}}

\newcommand{\Ker}{\operatorname{Ker}}
\newcommand{\Image}{\operatorname{Im}}

\newcommand{\ord}{\operatorname{ord}}
\newcommand{\Pic}{\operatorname{Pic}}

\newcommand{\Hom}{\operatorname{Hom}}

\newcommand{\outdeg}{{\operatorname{outdeg}}}

\begin{document}

\title[Graphs and Riemann surfaces]{Riemann-Roch and Abel-Jacobi theory on a finite graph}


\subjclass[2000]{05C38, 14H55.}

\author{Matthew Baker}
\address{School of Mathematics\\
Georgia Institute of Technology\\
Atlanta, Georgia 30332-0160\\
USA}
\email{mbaker@math.gatech.edu \\ snorine@math.gatech.edu}

\author{Serguei Norine}

\begin{abstract}
It is well-known that a finite graph can be viewed, in many respects, as a discrete analogue of
a Riemann surface.  In this paper, we pursue this analogy further in the context of linear equivalence
of divisors.  In particular, we formulate and prove a graph-theoretic analogue of the classical Riemann-Roch theorem.
We also prove several results, analogous to classical facts about Riemann surfaces,
concerning the Abel-Jacobi map from a graph to its Jacobian.
As an application of our results, we characterize the existence or non-existence of a winning
strategy for a certain chip-firing game played on the vertices of a graph.
\end{abstract}

\thanks{We would like to thank Robin Thomas for a number of useful discussions.
The first author would also like to thank his Summer 2006 REU student Dragos Ilas for computing
a number of examples and testing out conjectures about the graph-theoretic Abel-Jacobi map.
Thanks also to Hendrik Lenstra, Dino Lorenzini, and the anonymous referees for their helpful comments.
The first author's work was supported in part by NSF grant DMS-0600027, and the
second author's by NSF grant DMS-0200595.}

\maketitle

\section{Introduction}
\label{IntroSection}

\subsection{Overview}
\label{OverviewSection}

In this paper, we explore some new analogies between finite graphs and Riemann surfaces.
Our main result is a graph-theoretic analogue of the classical Riemann-Roch theorem.
We also study the Abel-Jacobi map $S$ from a graph $G$ to its Jacobian, as well as the
higher symmetric powers $S^{(k)}$ of $S$.
We prove, for example, that $S^{(g)}$ is always surjective, and that
$S^{(1)}$ is injective when $G$ is $2$-edge-connected.
These results closely mirror classical facts about the Jacobian of a Riemann surface.
As an application of our results, we characterize the existence or non-existence of a winning
strategy for a certain chip-firing game played on the vertices of a graph.

The paper is structured as follows.  In \S\ref{IntroSection}, we provide all of the relevant definitions and
state our main results.
The proof of the Riemann-Roch theorem for graphs occupies \S\ref{CombinatorialRiemannRochSection}-\ref{GraphRiemannRochSection}.
In \S\ref{AbelJacobiSection}, we study the injectivity and surjectivity of $S^{(k)}$ for $k\geq 1$, and
explain the connection with the chip-firing game.
Related results and further questions are discussed in \S\ref{ComplementsSection}.
The paper concludes with two appendices.
In Appendix~\ref{RSAppendix}, we provide the reader with a brief summary of some classical results about
Riemann surfaces, and in Appendix~\ref{GraphAbelTheoremSection}, we discuss the
graph-theoretic analogue of Abel's theorem proved in \cite{BDN}.

 \subsection{Notation and Terminology}
 \label{NotationSection}

Throughout this paper, a {\em Riemann surface} will mean a compact, connected one-dimensional complex manifold, and
a {\em graph} will mean a finite, unweighted multigraph having no loop edges.
All graphs in this paper are assumed to be connected.
We denote by $V(G)$ and $E(G)$, respectively, the set of vertices and edges of $G$.
We will simply write $G$ instead of $V(G)$ when there is no danger of confusion.
Also, we write $E_v = E_v(G)$ for the set of edges incident to a given vertex $v$.

For $k \geq 2$, a graph $G$ is called {\em $k$-edge-connected} if
$G-W$ is connected for every set $W$ of at most $k-1$ edges of $G$.
(By convention, we consider the trivial graph having one vertex and no edges to be
$k$-edge-connected for all $k$.)
Alternatively, define a {\em cut} to be the set of all edges connecting a vertex in $V_1$ to a vertex in $V_2$
for some partition
of $V(G)$ into disjoint non-empty subsets $V_1$ and $V_2$.
Then $G$ is $k$-edge-connected if and only if every cut has size at least $k$.



If $A \subseteq V(G)$, we denote by $\chi_A : V(G) \to \{ 0,1 \}$ the characteristic function
of $A$.

 \subsection{The Jacobian of a finite graph}
 \label{JacobianDefinitionSubsection}

Let $G$ be a graph, and choose an ordering $\{ v_1,\ldots, v_n \}$ of the vertices of $G$.
The {\em Laplacian matrix} $Q$ associated to $G$ is the $n \times n$
matrix $Q = D - A$, where $D$ is the diagonal matrix whose
$(i,i)^{\rm th}$ entry is the degree of vertex $v_i$, and $A$ is the {\em adjacency matrix} of the
graph, whose $(i,j)^{\rm th}$ entry is the number of edges connecting $v_i$ and $v_j$.
Since loop edges are not allowed, the $(i,i)^{\rm th}$ entry of $A$ is zero for all $i$.
It is well-known and easy to verify that $Q$ is symmetric, has rank $n-1$, and that the
kernel of $Q$ is spanned by the vector whose entries are all equal to 1 (see
\cite{BiggsBook,Bollobas,GR}).

Let $\Div(G)$ be the free abelian group on the set of vertices of $G$.
We think of elements of $\Div(G)$ as formal integer linear combinations of elements of $V(G)$,
and write an element $D \in \Div(G)$ as $\sum_{v \in V(G)} a_v (v)$, where each $a_v$ is an integer.
By analogy with the Riemann surface case, elements of $\Div(G)$ are called {\em divisors} on $G$.


\medskip

For convenience, we will write $D(v)$ for the coefficient $a_v$ of $(v)$ in $D$.

\medskip

There is a natural partial order on the group $\Div(G)$: we say that $D \geq D'$ if and only if
$D(v) \geq D'(v)$ for all $v \in V(G)$.  A divisor $E \in \Div(G)$ is called {\em effective} if
$E \geq 0$.  We write $\Div_+(G)$ for the set of all effective divisors on $G$.

The {\em degree function} $\deg : \Div(G) \to \ZZ$ is defined by
$\deg(D) = \sum_{v \in V(G)} D(v)$.

\begin{remark}
Note that the definitions of the partial order $\geq$, the space $\Div_+(G)$, and the map $\deg$
make sense when $V(G)$ is replaced by an arbitrary set $X$.  This observation will be used
in \S\ref{CombinatorialRiemannRochSection} when we formulate an abstract ``Riemann-Roch Criterion''.
\end{remark}

\medskip

We let $\M(G) = \Hom(V(G), \ZZ)$ be the abelian group consisting of all integer-valued functions on the vertices
of $G$.  One can think of $\M(G)$ as analogous to the field $\M(X)$ of meromorphic functions on a Riemann surface $X$
(though it is actually more like the abelian group $\{ \log |f| \; : \; f \in \M(X)^* \}$, see Remark~\ref{LaplacianDivisorRemark}).

Using our ordering of the vertices, we obtain isomorphisms between $\Div(G), \M(G)$, and the space of
$n \times 1$ column vectors having integer coordinates.
%
We write $[D]$ (resp. $[f]$) for the column vector corresponding to $D \in \Div(G)$ (resp. $f \in \M(G)$).
The {\em Laplacian operator} $\Delta : \M(G) \to \Div(G)$ is given by the formula
\[
\Delta(f) = \sum_{v \in V(G)} \Delta_v(f) (v) \ ,
\]
where
\[
\begin{aligned}
\Delta_v(f) &= \deg(v) f(v) - \sum_{e = wv \in E_v} f(w) \\
&= \sum_{e = wv \in E_v} (f(v) - f(w)) \ . \\
\end{aligned}
\]

In terms of matrices, it follows from the definitions that
\[
[\Delta(f)] = Q [f] \ .
\]

\begin{remark}
The fact that $Q$ is a symmetric matrix is equivalent to the fact that $\Delta$ is self-adjoint with respect
to the bilinear pairing $\la f, D \ra = \sum_{v \in V(G)} f(v) D(v)$ on $\M(G) \times \Div(G)$.
This is the graph-theoretic analogue of the {\em Weil reciprocity theorem} on a Riemann surface
(see p. 242 of \cite{GH} and Remark~\ref{LaplacianDivisorRemark} below).
\end{remark}

We define the subgroup $\Div^0(G)$ of $\Div(G)$ consisting of {\em divisors of degree zero} to be the kernel of $\deg$, i.e.,
\[
\Div^0(G) = \{ D \in \Div(G) \; : \; \deg(D) = 0 \} \ .
\]
More generally, for each $k \in \ZZ$ we define
$\Div^k(G) = \{ D \in \Div(G) \; : \; \deg(D) = k \}$, and
$\Div_+^k(G) = \{ D \in \Div(G) \; : \; D \geq 0 \textrm{ and } \deg(D) = k \}$.
The set $\Div_+^1(G)$ is canonically isomorphic to $V(G)$.

We also define the subgroup $\Prin(G)$ of $\Div(G)$ consisting of {\em principal divisors} to be
the image of $\M(G)$ under the Laplacian operator, i.e.,
\begin{equation}
\label{eq:PrinDef}
\Prin(G) := \Delta(\M(G)) \ .
\end{equation}

It is easy to see that every principal divisor has degree zero, so that $\Prin(G)$ is a subgroup of $\Div^0(G)$.

\begin{remark}
\label{LaplacianDivisorRemark}
The classical motivation for (\ref{eq:PrinDef}) is that the divisor of a nonzero meromorphic function
$f$ on a Riemann surface $X$ can be recovered
from the extended real-valued function $\log |f|$ using the (distributional) Laplacian operator $\Delta$.
More precisely if $\Delta(\varphi)$ is defined so that
\[
\int_X \psi \, \Delta(\varphi) = \int_X \varphi \, \Delta(\psi)
\]
for all suitably smooth test functions $\psi : X \to \RR$,
where $\Delta(\psi)$ is given in local coordinates by the formula
\[
\Delta(\psi) = \frac{1}{2\pi} \left( \frac{\partial^2 \psi}{\partial x^2}
+ \frac{\partial^2 \psi}{\partial y^2} \right) dx \wedge dy \ ,
\]
then
\[
\Delta(\log|f|) = \sum_{P \in X} \ord_P(f) \delta_P \ .
\]
In other words, the divisor of $f$ can be identified with the Laplacian of $\log |f|$.
\end{remark}

Following \cite{BDN} and \cite{Nagnibeda}, we define the group $\Jac(G)$, called the {\em Jacobian} of $G$, to be the corresponding
quotient group:
\begin{equation}
\label{JacobianDef}
\Jac(G) = \frac{\Div^0(G)}{\Prin(G)} \ .
\end{equation}
As shown in \cite{BDN}, $\Jac(G)$ is a finite abelian group whose order $\kappa(G)$ is the
number of {\em spanning trees} in $G$.  (This is a direct consequence of Kirchhoff's famous
{\em Matrix-Tree Theorem}, see \S{14} of \cite{BiggsAPTG}.)  The group $\Jac(G)$ is a
discrete analogue of the Jacobian of a Riemann surface.
We will write $[D]$ for the class in $\Jac(G)$ of a divisor $D \in \Div^0(G)$.  (There should not be any confusion
between this notation and our similar notation for the column vector associated to a divisor.)

\medskip

In \cite{BDN}, the group $\Jac(G)$ is called the {\em Picard group}, and denoted $\Pic(G)$, and the term {\em Jacobian} is
reserved for an {\em a priori} different group denoted $J(G)$.  However, as shown in
Proposition 7 of \cite{BDN}, the two groups are canonically isomorphic.
The isomorphism $\Pic(G) \cong J(G)$ is the graph-theoretic analogue of {\em Abel's theorem}
(see Theorem VIII.2.2 of \cite{Miranda}).


\subsection{The Abel-Jacobi map from a graph to its Jacobian}
\label{AbelJacobiSubsection}

If we fix a base point $v_0 \in V(G)$, we can define the {\em Abel-Jacobi map} $S_{v_0} : G \to \Jac(G)$
by the formula
\begin{equation}
\label{eq:AbelJacobiMap}
S_{v_0}(v) = [(v) - (v_0)] \ .
\end{equation}

We also define, for each natural number $k \geq 1$, a map $S_{v_0}^{(k)} : \Div_+^k(G) \to \Jac(G)$
by
\[
S_{v_0}^{(k)}((v_1) + \cdots + (v_k)) = S_{v_0}(v_1) + S_{v_0}(v_2) + \cdots + S_{v_0}(v_k) \ .
\]

The map $S_{v_0}$ can be characterized by the following universal property (see \S{3} of \cite{BDN}).
A map $\varphi : G \to A$ from $V(G)$ to an abelian group $A$ is called {\em harmonic} if
for each $v \in G$, we have
\[
\deg(v) \cdot \varphi(v) = \sum_{e=wv \in E_v} \varphi(w) \ .
\]
Then $S_{v_0}$ is universal among all harmonic maps from $G$ to abelian groups sending $v_0$ to $0$,
in the sense that if $\varphi : G \to A$ is any such map, then there is a unique group homomorphism
$\psi : \Jac(G) \to A$ such that $\varphi = \psi \circ S_{v_0}$.

\medskip

Let $g = |E(G)| - |V(G)| + 1$ be the {\em genus}\footnote{In graph theory, the term ``genus'' is traditionally used for a different concept,
namely, the smallest genus (i.e., first Betti number) of any surface in which the graph can be embedded, and the integer $g$ is called the ``cyclomatic number'' of $G$.
We call $g$ the genus of $G$ in order to highlight the analogy with Riemann surfaces.}
of $G$, which is the number of linearly independent cycles of $G$, or equivalently, the first Betti number of $G$ (i.e., the dimension of
$H_1(G, \RR)$).

We write $S$ instead of $S_{v_0}$ when the base point $v_0$ is understood.
In \S\ref{AbelJacobiSection}, we will prove:

\begin{theorem}
\label{SurjectivityTheorem}
The map $S^{(k)}$ is surjective if and only if $k \geq g$.
\end{theorem}

The surjectivity of $S^{(g)}$ is the graph-theoretic analogue of a classical result about Riemann surfaces known as Jacobi's Inversion
Theorem (see p.~235 of \cite{GH}).  For a Riemann surface $X$, it is clear that $S^{(g-1)} : X^{(g-1)} \to \Jac(X)$ is not
surjective, since $\dim S^{(g-1)} = g-1 < \dim \Jac(X) = g$.

\medskip

As a complement to Theorem~\ref{SurjectivityTheorem}, we will also precisely characterize the values of $k$ for which
$S^{(k)}$ is injective:

\begin{theorem}
\label{InjectivityTheorem}
The map $S^{(k)}$ is injective if and only if $G$ is $(k+1)$-edge-connected
\end{theorem}

For $2$-edge-connected graphs, Theorem~\ref{InjectivityTheorem} is the analogue of the well-known fact
that the Abel-Jacobi map from a Riemann surface $X$ to its Jacobian is injective if and only if $X$ has
genus at least 1.
(See Proposition VIII.5.1 of \cite{Miranda}.)



 \subsection{Chip-firing games on graphs}
 \label{ChipFiringSubsection}

There have been a number of papers devoted to ``chip-firing games'' played on the vertices of a graph;
see, e.g., \cite{BiggsCF,BL,BLS,GR,Merino1,Merino2,Tardos,vdH}.
In this paper, as an application of Theorem~\ref{SurjectivityTheorem},
we study a new chip firing game with some rather striking features.

Our chip-firing game, like the one considered by Biggs in \cite{BiggsCF} (see also \S{31-32} of \cite{BiggsAPTG}),
is most conveniently stated using ``dollars'' rather than chips.
Let $G$ be a graph, and consider the following game of
``solitaire'' played on the vertices of $G$.  The initial configuration of the game assigns to each vertex $v$ of $G$
an integer number of dollars.  Such a configuration can be identified with a divisor $D \in \Div(G)$.
A vertex which has a negative number of dollars assigned to it is said to be in {\em debt}.
A {\em move} consists of a vertex $v$ either borrowing one dollar from each of its neighbors or giving one dollar
to each of its neighbors.  Note that any move leaves the total number of dollars unchanged.
The object of the game is to reach, through a sequence of moves,
a configuration in which no vertex is in debt.  We will call such a configuration a {\em winning position},
and a sequence of moves which achieves such a configuration a {\em winning strategy}.

As before, we let $g = |E(G)| - |V(G)| + 1$.
In \S\ref{ChipFiringSection}, we will prove
the following result by showing that it is equivalent to Theorem~\ref{SurjectivityTheorem}:

\begin{theorem}
\label{ChipFiringTheorem}
Let $N = \deg(D)$ be the total number of dollars present at any stage of the game.
\begin{itemize}
\item[1.] If $N \geq g$, then there is always a winning strategy.
\item[2.] If $N \leq g-1$, then there is always an initial configuration for which
no winning strategy exists.
\end{itemize}
\end{theorem}

See \S\ref{BLSSection} for a discussion of the relationship between our chip-firing game and the one studied by
Bj{\"o}rner, Lov{\'a}sz, and Shor in \cite{BLS}, and see \S\ref{BiggsSection} for a discussion of the relationship between
our chip-firing game and the dollar game of Biggs.

 \subsection{Linear systems and the Riemann-Roch theorem}

We define an equivalence relation $\sim$ on the group $\Div(G)$ by declaring that $D \sim D'$ if and only if
$D - D' \in \Prin(G)$.  Borrowing again from the theory of Riemann surfaces,
we call this relation {\em linear equivalence}.
Since a principal divisor has degree zero, it follows that linearly equivalent divisors have the same degree.
Note that by (\ref{JacobianDef}), the Jacobian of $G$ is the set of
linear equivalence classes of degree zero divisors on $G$.

For $D \in \Div(G)$, we define the {\em linear system associated to $D$} to
be the set $|D|$ of all effective divisors linearly equivalent to $D$:
\[
|D| = \{ E \in \Div(G) \; : \; E \geq 0, \; E \sim D \} \ .
\]

As we will see in \S\ref{ChipFiringSection}, it follows from the definitions that two divisors $D$ and $D'$ on $G$ are linearly equivalent
if and only if there is a sequence of moves taking $D$ to $D'$
in the chip firing game described in \S\ref{ChipFiringSubsection}.  It follows that there is a winning strategy in the
chip-firing game whose initial configuration corresponds to $D$ if and only if $|D| \neq \emptyset$.

\medskip

We define the {\em dimension} $r(D)$ of the linear system $|D|$ by setting $r(D)$ equal to $-1$ if
$|D| = \emptyset$, and then declaring that for each integer $s \geq 0$,
$r(D) \geq s \textrm{ if and only if } |D - E| \neq \emptyset$
for all effective divisors $E$ of degree $s$.
It is clear that $r(D)$ depends only on the linear equivalence class of $D$.

\begin{remark}
By Lemma~\ref{ChipFiringEquivalenceLemma} below, we have $r(D) \geq 0$ if and only if there is a winning strategy in the
chip firing game with initial configuration $D$, $r(D) \geq 1$ if and only if there is still a winning strategy
after subtracting one dollar from any vertex, etc.
\end{remark}

The {\em canonical divisor} on $G$ is the divisor $K$ given by
\begin{equation}
\label{eq:canondiv}
K = \sum_{v \in V(G)} \left( \deg(v) - 2 \right) (v) \ .
\end{equation}

Since the sum over all vertices $v$ of $\deg(v)$ equals twice the number of edges in $G$,
we have $\deg(K) = 2|E(G)| - 2|V(G)| = 2g - 2$.

\medskip

We can now state a graph-theoretic analogue of the Riemann-Roch theorem
(see Theorem VI.3.11 of \cite{Miranda}).
The proof will be given in \S\ref{GraphRiemannRochSection}.

\begin{theorem}[Riemann-Roch for Graphs]
\label{GraphRiemannRochTheorem}
Let $G$ be a graph, and let $D$ be a divisor on $G$.  Then
\[
r(D) - r(K - D) = \deg(D) + 1 - g \ .
\]
\end{theorem}

\begin{remark}
\label{DimensionRemark}

(i) Our definition of $r(D)$ agrees with the usual definition of $r(D)$ as $\dim L(D) - 1$ in the Riemann surface case
(see, e.g., p.~250 of \cite{GH} or \S{III.8.15} of \cite{FarkasKra}).

(ii) One must be careful, however, not to rely too much on intuition from the Riemann surface case
when thinking about the quantity $r(D)$ for divisors on graphs.
For example, for Riemann surfaces one has $r(D)=0$ if and only if $|D|$ contains exactly one element,
but neither implication is true in general for graphs.  For example, consider the
canonical divisor $K$ on a graph $G$ with two vertices $v_1$ and
$v_2$ connected by $m$ edges. Then clearly $r(K) \geq m-2$, and in
fact we have $r(K) = m-2$.  (This can be proved directly, or deduced as a
consequence of Theorem~\ref{GraphRiemannRochTheorem}.) However, $|K| =
\{K\}$ as
\[
D \sim K \, \Leftrightarrow \, \exists \, i \in \ZZ \; : \;
D = (m-2+im)(v_1)+(m-2-im)(v_2) \ .
\]

To see that the other implication also fails,
consider a graph $G$ with $V(G) = \{ v_1, v_2, v_3, v_4, v_5 \}$,
$E(G)=\{v_1v_2, v_2v_3,v_3v_4,v_4v_5,v_5v_1,v_3v_1\}$, and $D =
2(v_4)\in \Div(G)$. Then $(v_3)+(v_5) \in |D|$, but it follows from
Lemma ~\ref{OrderDivisorsLemma} (or can be verified directly) that $|D-(v_1)| = \emptyset$,
and therefore $r(D)=0$.


(iii) The set $L(D) := \{ f \in \M(G) \; : \; \Delta(f) \geq -D \}$ is not a vector space, so one cannot just define the
number $r(D)$ as $\dim L(D) - 1$ as in the classical case.  This should not be surprising, since elements of $L(D)$ are
analogous to functions of the form $\log |f|$ with $f$ a nonzero meromorphic function on a Riemann surface $X$.  On the other hand,
$L(D) \cup \{ \infty \}$ is naturally a finitely generated semimodule over
the tropical semiring $(\NN \cup \{ \infty \},\min,+)$ (see \S2.4 of \cite{Gaubert}), and there is a natural notion in this context
for the {\em dimension} of $L(D)$ (see Corollary 95 in \cite{Gaubert}).
However, examples like the ones above show that the tropical dimension of $L(D)$ is {\em not} the same as $r(D) + 1$, and does not
obey Theorem~\ref{GraphRiemannRochTheorem}.

\end{remark}

\section{A Riemann-Roch criterion}
\label{CombinatorialRiemannRochSection}

In this section, we formulate an abstract criterion giving necessary and sufficient conditions for the
Riemann-Roch formula $r(D) - r(K - D) = \deg(D) + 1 - g$
to hold, where $r(D)$ is defined in terms of an equivalence relation
on an arbitrary free abelian group.  This result, which is purely combinatorial in nature,
will be used in \S\ref{GraphRiemannRochSection} in our proof of the Riemann-Roch theorem for graphs.

\medskip

The general setup for our result is as follows.

\medskip

Let $X$ be a non-empty set, and let $\Div(X)$ be the free abelian group on $X$.
As usual, elements of $\Div(X)$ are called {\em divisors} on $X$,
divisors $E$ with $E \geq 0$ are called {\em effective},
and for each integer $d$, we denote by $\Div_+^{d}(X)$ the set of effective divisors of degree $d$ on $X$.

Let $\sim$ be an equivalence relation on $\Div(X)$ satisfying the following two properties:
\begin{itemize}
\item[(E1)] If $D \sim D'$ then $\deg(D) = \deg(D')$.
\item[(E2)] If $D_1 \sim D_1'$ and $D_2 \sim D_2'$, then $D_1 + D_1' \sim D_2 + D_2'$.
\end{itemize}

For each $D \in \Div(X)$, define $|D| = \{ E \in \Div(G) \; : \; E \geq 0, \; E \sim D \}$,
and define the function $r : \Div(X) \to \{ -1, 0, 1, 2, \ldots \}$ by declaring that
for each integer $s \geq 0$,
\[
r(D) \geq s \iff |D - E| \neq \emptyset \, \forall \, E \in \Div(X)
\; : \; E \geq 0 \textrm{ and } \deg(E) = s \ .
\]
Note that the above equivalence is true for all integers $s$.
It is easy to see that $r(D) = -1$ if $\deg(D) < 0$, and if $\deg(D) = 0$ then
$r(D) = 0$ if $D \sim 0$ and $r(D) = -1$ otherwise.

\begin{lemma}
\label{SubAdditivityLemma}
For all $D,D' \in \Div(X)$ such that
$r(D),r(D') \geq 0$, we have $r(D + D') \geq r(D) + r(D')$.
\end{lemma}

\begin{proof}
Let $E_0 = (x_1) + \cdots + (x_{r(D)+r(D')})$ be an arbitrary effective divisor of
degree $r(D) + r(D')$, and let
$E = (x_1) + \cdots + (x_{r(D)})$ and
$E' = (x_{r(D)+1}) + \cdots + (x_{r(D) + r(D')})$.
Then $|D - E|$ and $|D' - E'|$ are non-empty, so that
$D - E \sim F$ and $D' - E' \sim F'$ with $F,F' \geq 0$.
It follows that $(D + D') - (E + E') = (D + D') - E_0 \sim F + F' \geq 0$,
and thus $r(D+D') \geq r(D) + r(D')$.
\end{proof}

Let $g$ be a nonnegative integer, and define
\[
\N = \{ D \in \Div(X) \; : \; \deg(D) = g-1 \textrm{ and } |D| = \emptyset \} \ .
\]

Finally, let $K$ be an element of $\Div(X)$ having degree $2g-2$.
The following theorem gives necessary and sufficient conditions for the Riemann-Roch formula
to hold for elements of $\Div(X) / \sim$.

\begin{theorem}
\label{CombinatorialRiemannRochTheorem}
Define $\epsilon : \Div(X) \to \ZZ / 2\ZZ$ by declaring that $\epsilon(D) = 0$ if
$|D| \neq \emptyset$ and $\epsilon(D) = 1$ if $|D| = \emptyset$.
Then the Riemann-Roch formula
\begin{equation}
\label{RRFormula}
r(D) - r(K - D) = \deg(D) + 1 - g
\end{equation}
holds for all $D \in \Div(X)$ if and only if the following two
properties are satisfied:
\begin{itemize}
\item[(RR1)] For every $D \in \Div(X)$, there exists $\nu \in \N$ such that
$\epsilon(D) + \epsilon(\nu - D) = 1$.
\item[(RR2)] For every $D \in \Div(X)$ with $\deg(D) = g-1$, we have
$\epsilon(D) + \epsilon(K-D) = 0$.
\end{itemize}
\end{theorem}

\begin{remark}
\label{RRRemark}

(i) Property (RR2) is equivalent to the assertion that $r(K) \geq g-1$.
Indeed, if (RR2) holds then for every effective divisor $E$ of degree $g-1$,
we have $|K - E| \neq \emptyset$, which means that $r(K) \geq g-1$.
Conversely, if $r(K) \geq g-1$ then $\epsilon(K-E) = \epsilon(E) = 0$
for every effective divisor $E$ of degree $g-1$.  Therefore
$\epsilon(D) = 0$ implies $\epsilon(K-D) = 0$.  By symmetry,
we obtain $\epsilon(D) = 0$ if and only if $\epsilon(K-D) = 0$, which is
equivalent to (RR2).

(ii) When the Riemann-Roch formula (\ref{RRFormula}) holds, we automatically have $r(K) = g-1$.

\end{remark}

\begin{remark}

(i) When $X$ is a Riemann surface and $\sim$ denotes linear equivalence of divisors,
then one can show independently of the Riemann-Roch theorem that $r(K) = g-1$, i.e., that
the vector space of holomorphic $1$-forms on $X$ is $g$-dimensional.  Thus one can prove
directly that (RR2) holds.  We do not know if there is a direct proof of (RR1)
which does not make use of Riemann-Roch, but if so, one could deduce the classical Riemann-Roch
theorem from it using Theorem~\ref{CombinatorialRiemannRochTheorem}.

(ii) Divisors of degree $g-1$ on a Riemann surface $X$ which belong to $\N$ are classically referred to as {\em non-special}
(which explains our use of the symbol $\N$).
\end{remark}

Before giving the proof of Theorem~\ref{CombinatorialRiemannRochTheorem}, we need a couple of
preliminary results.  The first is the following simple lemma, whose proof is left to the reader.

\begin{lemma}
\label{MinimizationLemma}
Suppose $\psi : A \to A'$ is a bijection between sets, and that $f : A \to \ZZ$ and
$f' : A' \to \ZZ$ are functions which are bounded below.
If there exists a constant $c \in \ZZ$ such that
\[
f(a) - f'(\psi(a)) = c
\]
for all $a \in A$, then
\[
\min_{a \in A} f(a) - \min_{a' \in A'} f'(a') = c \ .
\]
\end{lemma}

\medskip

If $D = \sum_i a_i (x_i) \in \Div(X)$, we define
\[
\deg^+(D) = \sum_{a_i \geq 0} a_i \ .
\]

The key observation needed to deduce (\ref{RRFormula}) from (RR1) and (RR2) is
the following alternate characterization of the quantity $r(D)$:

\begin{lemma}
\label{AltDimensionLemma}
If (RR1) holds then for every $D \in \Div(X)$ we have
\begin{equation}
\label{DimensionFormula}
r(D) = \left( \min_{\substack{D' \sim D \\ \nu \in \N} } \deg^+(D' - \nu) \right) - 1 \ .
\end{equation}
\end{lemma}

\begin{proof}
Let $r'(D)$ denote the
right-hand side of (\ref{DimensionFormula}).
If $r(D) < r'(D)$, then there exists an effective divisor $E$ of degree $r'(D)$ for which
$r(D - E) = -1$.
By (RR1), this means that there exists a divisor $\nu \in \N$ and an effective divisor $E'$ such that
$\nu - D + E \sim E'$.  But then $D' - \nu = E - E'$ for some divisor $D' \sim D$,
and thus
\[
\deg^+(D' - \nu) - 1 \leq \deg(E) - 1 = r'(D) - 1 \ ,
\]
contradicting the definition of $r'(D)$.  It follows that $r(D) \geq r'(D)$.

Conversely, if we choose divisors $D' \sim D$ and $\nu \in \N$ achieving the minimum
in (\ref{DimensionFormula}), then $\deg^+(D' - \nu) = r'(D) + 1$, and therefore
there are effective divisors $E,E'$ with $\deg(E) = r'(D) + 1$ such that
$D' - \nu = E - E'$.  But then $D - E \sim \nu - E'$, and since $\nu - E'$ is not
equivalent to any effective divisor, it follows that $|D - E| = \emptyset$.  Therefore
$r(D) \leq r'(D)$.
\end{proof}

We can now give the proof of Theorem~\ref{CombinatorialRiemannRochTheorem}.

\begin{proof}[Proof of Theorem~\ref{CombinatorialRiemannRochTheorem}]

We first prove that (\ref{RRFormula}) implies (RR1) and (RR2).

Let $D$ be a divisor on $X$, and let $d=\deg(D)$.
Property (RR2) is more or less immediate, since (\ref{RRFormula}) implies that
if $\deg(D) = g-1$ then $r(D) = r(K - D)$.

We cannot have $\epsilon(D) = \epsilon(\nu - D) = 0$,
or else by Lemma~\ref{SubAdditivityLemma} we would have
$r(\nu) \geq 0$, contradicting the definition of $\N$.
As we will see in the next paragraph, $\N$ is non-empty; 
therefore, to prove (RR1) it suffices to show that if $r(D) = -1$ then
$r(\nu - D) \geq 0$ for some $\nu \in \N$.

If $r(D+E) \geq 0$ for all $E \in \Div_+^{g-1-d}(X)$, then
(\ref{RRFormula}) implies that $r(K-D-E) \geq 0$ for all such $E$,
and therefore $r(K-D) \geq g-1-d$.  Another application of
(\ref{RRFormula}) then yields $r(D) = r(K-D) + d + 1 - g \geq 0$.

Therefore when $r(D) = -1$, there exists an effective divisor $E$ of degree $g-1-d$
such that $r(D+E) = -1$.  Since $\deg(D+E) = g-1$, this means that $D+E \in \N$, and
therefore $D + E = \nu$ for some $\nu \in \N$.  For this choice of $\nu$, we have
$r(\nu-D) \geq 0$, which proves (RR1).

\medskip

We now show that (RR1) and (RR2) imply (\ref{RRFormula}).
Let $D \in \Div(X)$.  For every $\nu \in \N$, property (RR2) implies that
$\overline{\nu} := K - \nu$ is also in $\N$.
Writing
\[
\nu - D' = K - D' - \overline{\nu} \ ,
\]
it follows that
\begin{equation*}
\label{VariationalEquation}
\begin{aligned}
\deg^+(D' - \nu) - \deg^+((K - D') - \overline{\nu})
     &= \deg^+(D' - \nu) - \deg^+(\nu - D') \\
     &= \deg(D' - \nu) \\
     &= \deg(D) + 1 - g \ . \\
\end{aligned}
\end{equation*}

Since the difference $\deg^+(D' - \nu) - \deg^+((K - D') - \overline{\nu})$ has the constant value
$\deg(D) + 1 - g$ for all $D'$ and $\nu$,
and since $\overline{\nu} = K - \nu$ runs through all possible elements of $\N$ as $\nu$ does,
it follows from Lemmas \ref{MinimizationLemma} and \ref{AltDimensionLemma} that
$r(D) - r(K - D) = \deg(D) + 1 - g$ as desired.
\end{proof}

\section{Riemann-Roch for graphs}
\label{GraphRiemannRochSection}

 \subsection{$G$-parking functions and reduced divisors}
 \label{GParkingSection}

In this section, we use the notion of a {\em $G$-parking function},
introduced in ~\cite{PS}, to define a unique {\em reduced divisor} in each
equivalence class in $\Div(G)$. Reduced divisors will play a key role
in our proof of the Riemann-Roch theorem for graphs in the next section.
Our reduced divisors are closely related to the ``critical configurations''
considered by Biggs in~\cite{BiggsCF,BiggsTutte}, as we will explain in Section~\ref{BiggsSection}.

We now present the relevant definitions. For $A \subseteq V(G)$ and
$v \in A$, let $\outdeg_A(v)$ denote the number of edges of $G$
having $v$ as one endpoint and whose other endpoint lies in $V(G) - A$.
Select a vertex $v_{0} \in V(G)$. We say that a function $f:V(G)-
\{v_0\} \longrightarrow \ZZ$ is a \emph{$G$-parking function}
(relative to the base vertex $v_0$) if the following two conditions are
satisfied:
\begin{itemize}
\item[(P1)] $f(v) \geq 0$ for all $v \in V(G)- \{v_0\}$.
\item[(P2)] For every non-empty set $A \subseteq V(G)- \{v_0\}$, there exists
a vertex $v \in A$ such that $f(v) < \outdeg_A(v)$.
\end{itemize}

We say that a divisor $D \in \Div(G)$ is \emph{$v_0$-reduced} if the map $v \mapsto D(v)$, defined
for $v \in V(G)- \{v_0\}$, is a $G$-parking function.
In terms of the chip-firing game, a divisor $D$ is $v_0$-reduced if and only if (1) no vertex $v \neq v_0$ is in debt; and 
(2) for every non-empty subset $A$ of $V(G)- \{v_0\}$, if all vertices in $A$ were to perform a lending move, some vertex in $A$ would go into debt.

\begin{proposition}
\label{GParkingProp}
Fix a base vertex $v_0 \in V(G)$.  Then
for every $D \in \Div(G)$, there exists a unique $v_0$-reduced divisor $D' \in
\Div(G)$ such that $D' \sim D$.
\end{proposition}

\begin{proof}
We begin by presenting an informal sketch of the proof that such a divisor $D'$
exists in terms of the chip-firing game. We need to show that any
initial configuration can be transformed into a configuration corresponding
to a $v_0$-reduced divisor via a sequence of legal moves.
To accomplish this, we first obtain a configuration where no vertex except $v_0$ is in debt.
This can be done, for example, by arranging the vertices in some order,
starting with $v_0$, in such a way that every vertex except for
$v_0$ has a neighbor that precedes it in this order. We then take
the vertices out of debt consecutively, starting with the last
vertex, by at each step having some neighbor $w$ which precedes the current vertex $v$
in the designated order lend out enough money to take $v$ out of debt.

Once we have obtained a configuration where no vertex other than $v_0$
is in debt, we enumerate the non-empty subsets $A_1,\ldots,A_s$ of $V(G)- \{v_0\}$.
If every vertex of $A_1$ can give a dollar to each of its neighbors
outside $A_1$ and remain out of debt, then each vertex of $A_1$ does so
(this is a combination of legal moves in the chip-firing
game); otherwise, we move on to the next set $A_2$, and so on.
Once the vertices in some set $A_i$ lend out money, we cycle through the
entire procedure again, beginning with $A_1$.
If for each $1\leq i \leq s$, there is some vertex in $A_i$ which cannot lend a dollar
to each of its neighbors outside $A_i$ without going into debt, then the
procedure terminates.

Note that $v_0$ never lends money during this procedure, and so it
must stop receiving money at some point.  None of the neighbors of $v_0$ lend
money out from this point on, and so they, too, must eventually stop
receiving money.  Iterating this argument, we see that the entire procedure
has to stop. The configuration $D'$ obtained at the end of this process
corresponds to a $v_0$-reduced divisor.

\medskip

We now formalize the argument presented above. For a vertex $v \in
V(G)$, let $d(v)$ denote the length of the shortest path in $G$
between $v$ and $v_0$. Let $d =  \max_{v \in V(G)}d(v)$ and let $S_k
= \{v \in V(G) \; : \; d(v)=k\}$ for $0 \leq k \leq d$.

Define the vectors $\mu_1(D)\in \ZZ^{d}$ and $\mu_2(D) \in
\ZZ^{d+1}$ by
\[
\mu_1(D) = \left( \sum_{\substack{v \in S_d \\ D(v)<0}} D(v),
\sum_{\substack{v \in S_{d-1} \\ D(v)<0}} D(v),\ldots,
\sum_{\substack{v \in S_{1} \\ D(v)<0}} D(v) \right) \ ,
\]
\[
\mu_2(D) = \left( \sum_{v \in S_0} D(v), \sum_{v \in S_{1}}
D(v),\ldots, \sum_{v \in S_d} D(v)\right) \ .
\]
Replacing $D$ by an equivalent divisor if necessary, we may assume without loss of generality that $$\mu_1(D)= \max_{D'
\sim D}\mu_1(D') \ \ \text{and} \ \ \mu_2(D)= \max_{\substack{D' \sim D \\
\mu_1(D) = \mu_1(D')}}\mu_2(D'),$$ where the maxima are taken in the
lexicographic order. It is easy to see that both maxima are
attained. We claim that the resulting divisor $D$ is $v_0$-reduced.

Suppose $D(v) < 0$ for some vertex $v \neq v_0$. Let $v'$ be a
neighbor of $v$ such that $d(v') < d(v)$ and let $D' = D -
\Delta(\chi_{\{v'\}})$. Then $D'(v)>D(v)$, and $D'(w) \geq D(w)$ for
every $w$ such that $d(w) \geq d(v)$. It follows that $\mu_1(D')
> \mu_1(D)$, contradicting the choice of $D$. Therefore $D(v) \geq
0$ for every $v \in V(G), v \neq v_0$.

Suppose now that for some non-empty subset $A \subseteq V(G)- \{v_0\}$, we
have $D(v) \geq \outdeg_A(v)$ for every $v \in A$. Let $D' = D -
\Delta(\chi_{A})$ and $d_A = \min_{v \in A}d(v)$. We have $D'(v)
\geq D(v)$ for all $v \in V(G) - A$ and $D'(v)=D(v)-\outdeg_A(v)
\geq 0$ for every $v \in A$. Therefore $\mu_1(D')=\mu_1(D)$, as they
are both the zero vector.  There must be a vertex $v' \in V(G)$ such that $d(v')<
d_A$, and for which $v'$ has a neighbor in $A$. It follows that $D'(v')>D(v')$,
and consequently $\mu_2(D') > \mu_2(D)$, once again contradicting
the choice of $D$. This finishes the proof of the claim.

It remains to show that distinct $v_0$-reduced divisors cannot be
equivalent. Suppose for the sake of contradiction that we are given
$v_0$-reduced divisors $D$ and $D'$ such that $D \sim D'$ and $D
\neq D'$. Let $f \in \M(G)$ be a function for which $D' - D =
\Delta(f)$. Then $f$ is non-constant, and by symmetry
we may assume that $f(v)>f(v_0)$ for some $v \in V(G)$. Let $A$
be the set of all the vertices $v \in V(G)$ for which $f(v)$ is
maximal. Then $v_0 \not\in A$, and for every $v \in A$ we have
$$
0 \leq D(v) =  D'(v) - \sum_{e = vw \in E_v} (f(v)-f(w)) \leq D'(v)
- \outdeg_A(v) \ .
$$
Thus $D'(v) \geq \outdeg_A(v)$ for every $v \in A$, contradicting
the assumption that $D'$ is $v_0$-reduced.
\end{proof}

\subsection{Proof of the Riemann-Roch theorem}
\label{GraphRiemannRochProofSection}

By Theorem~\ref{CombinatorialRiemannRochTheorem}, in order to prove
the Riemann-Roch theorem for graphs
(Theorem~\ref{GraphRiemannRochTheorem}), it suffices to verify
properties (RR1) and (RR2) when $X = G$ is a graph and $\sim$
denotes linear equivalence of divisors. This will be accomplished by
analyzing a certain family of divisors of degree $g-1$ on $G$.

For each linear (i.e., total) order $<_P$ on $V(G)$, we define
\[
\nu_{P} = \sum_{v \in V(G)}(|\{e = vw \in E(G) \; : \; w <_P v \}| - 1)(v).
\]

It is clear that $\deg(\nu_P)= |E(G)|-|V(G)|= g-1$.

\begin{lemma}
\label{OrderDivisorsLemma} For every linear order $<_P$ on $V(G)$ we
have $\nu_{P} \in \N$.
\end{lemma}
\begin{proof} Let $D \in \Div(G)$ be any divisor
of the form $D = \nu_{P} - \Delta(f)$ for some $f \in \M(G)$.
Let $V_f^{\max}$ be the set of vertices $v \in G$ at which $f$ achieves its
maximum value, and let $u$ be the minimal element of
$V_f^{\max}$ with respect to the order $<_P$.
Then $f(w) \leq f(u)$ for all $w \in V(G)$, and if $w <_P u$ then $f(w) < f(u)$.
Thus
\[
\begin{aligned}
D(u) &= \left( |\{ e = uw \in E(G) \; : \; w <_P u \}| - 1 \right)
- \sum_{e = uw \in E(G)} \left( f(u) - f(w) \right) \\
&= -1 + \sum_{\substack{e = uw \in E(G) \\ u <_P w} }(f(w) - f(u)) +
\sum_{\substack{e = uw \in E(G) \\ w <_P u}}(f(w) - f(u) + 1) \\
& \leq -1 \ ,
\end{aligned}
\]
since each term in these sums is non-positive by the choice of $u$.
It follows that $\nu_P$ is not equivalent to any effective divisor.
\end{proof}

\begin{theorem}
\label{MainGraphTheorem}
For every $D \in \Div(G)$, exactly one of
the following holds
\begin{itemize}
\item[(N1)] $r(D) \ge 0$; or
\item[(N2)] $r(\nu_P - D) \ge 0$ for some order $<_P$ on $V(G)$.
\end{itemize}
\end{theorem}

\begin{proof}
Choose $v_0 \in V(G)$. By Proposition~\ref{GParkingProp}, we may assume that
$D$ is $v_0$-reduced. We define $v_1,v_2, \ldots, v_{|V(G)|-1}$
inductively as follows. If $v_0,v_1, \ldots,v_{k-1}$ are defined,
let $A_k = V(G) - \{v_0,v_1, \ldots,v_{k-1}\}$, and let $v_k \in A_k$
be chosen so that $D(v_k) < \outdeg_{A_k}(v_k)$. Let $<_P$ be the
linear order on $V(G)$ such that $v_i <_P v_j$ if and only if $i<j$.

For every $1 \leq k \leq |V(G)|-1$ we have
\begin{align*}
D(v_k) &\leq \outdeg_{A_k}(v_k) - 1 \\ &= |\{e = v_kv_j \in E(G) \;
: \; j < k \}| -1 \\
&= \nu_P(v_k) \ .
\end{align*}
If $D(v_0) \geq 0$ then we have $D \ge 0$ and (N1) holds. If, on the
other hand, $D(v_0) \leq -1$ then $D \le \nu_P$ and (N2) holds.
Finally, note that if $r(D)\geq 0$ and $r(\nu_P - D) \geq 0$, then
$r(\nu_P) \geq 0$ by Lemma~\ref{SubAdditivityLemma}, contradicting
Lemma~\ref{OrderDivisorsLemma}.
\end{proof}

As a consequence of Lemma~\ref{OrderDivisorsLemma} and Theorem~\ref{MainGraphTheorem},
we obtain:

\begin{corollary}
\label{GraphNonSpecialCorollary}
 For $D \in \Div(G)$ with $\deg(D) = g -1$ we have $D \in \N$ if and only if
there exists a linear order $<_P$ on $V(G)$ such that $D \sim
\nu_P$.
\end{corollary}

\begin{proof}
It suffices to note that if $\nu_P - D \sim E$ with $E \geq 0$, then $\deg(E) = 0$
and thus $E = 0$, so that $D \sim \nu_P$.
\end{proof}

We can now prove our graph-theoretic version of the Riemann-Roch theorem.

\begin{proof}[Proof of Theorem~\ref{GraphRiemannRochTheorem}]
By Theorem~\ref{CombinatorialRiemannRochTheorem}, it suffices to
prove that conditions (RR1) and (RR2) are satisfied.

Let $D \in \Div(G)$, and suppose first that $r(D) \geq 0$.
Then for every $\nu \in \N$ we have $r(\nu - D) = -1$, and hence
$\epsilon(D) + \epsilon(\nu-D) = 0 + 1 = 1$ and (RR1) holds.
Suppose, on the other hand, that $r(D)<0$.
Then by Theorem~\ref{MainGraphTheorem}, we must have
$r(\nu_P - D) \ge 0$ for some order $<_P$ on $V(G)$, and then
$\epsilon(D) + \epsilon(\nu_P - D) = 1 + 0 = 1$.
As $\nu_P \in \N$ by Lemma ~\ref{OrderDivisorsLemma},
it follows once again that (RR1) holds.

To prove (RR2), it suffices to show that for every $D \in \N$ we have
$K-D \in \N$. By Corollary ~\ref{GraphNonSpecialCorollary}, we have
$D \sim \nu_P$ for some linear order $<_P$ on $V(G)$. Let $\bar P$
be the {\em reverse} of $P$ (i.e., $v <_P w \Leftrightarrow w <_{\bar P}
v$). Then for every $v \in V(G)$, we have
\begin{align*}
\nu_P(v) + \nu_{\bar P}(v)
= (|\{e = vw \in E(G) \; : \; w <_P v \}|
- 1) \\+ (|\{e = vw \in E(G) \; : \; w <_{\bar P} v \}| - 1) \\= \deg(v) -
2 = K(v) \ .
\end{align*}
 Therefore $K - D \sim K - \nu_P = \nu_{\bar P} \in \N$.
\end{proof}

 \subsection{Consequences of the Riemann-Roch theorem}

As in the Riemann surface case, one can derive a number of interesting consequences from
the Riemann-Roch formula.  As just one example, we prove a graph-theoretic analogue of
Clifford's theorem (see Theorem VII.1.13 of \cite{Miranda}).  For the statement, we call a divisor $D$
{\em special} if $|K-D| \neq \emptyset$, and {\em non-special} otherwise.

\begin{corollary}[Clifford's Theorem for Graphs]
\label{CliffordCor}
Let $D$ be an effective special divisor on a graph $G$.  Then
\[
r(D) \leq \frac{1}{2} \deg(D) \ .
\]
\end{corollary}

\begin{proof}
If $D$ is effective and special, then $K-D$ is also effective, and by Lemma~\ref{SubAdditivityLemma}
we have
\[
r(D) + r(K-D) \leq r(K) = g - 1 \ .
\]
On the other hand, by Riemann-Roch we have
\[
r(D) - r(K-D) = \deg(D) + 1 - g \ .
\]
Adding these two expressions gives $2 r(D) \leq \deg(D)$ as desired.
\end{proof}

As pointed out in \S{IV.5} of \cite{Hartshorne},
the interesting thing about Clifford's theorem is that for a non-special divisor $D$,
we can compute $r(D)$ exactly as a function of $\deg(D)$ using Riemann-Roch.  However,
for a special divisor, $r(D)$ does not depend only on the degree.  Therefore it is useful to
have a non-trivial upper bound on $r(D)$, and this is what Corollary~\ref{CliffordCor} provides.

\section{The Abel-Jacobi map from a graph to its Jacobian}
\label{AbelJacobiSection}

Let $G$ be a graph, let $v_0 \in V(G)$ be a base point, and let $k$ be a positive integer.
In this section, we discuss the injectivity and surjectivity of the
map $S_{v_0}^{(k)}$.

We leave it to the reader to verify the following elementary observations:

\begin{lemma}
\label{ElementaryObservationsLemma}
\begin{itemize}
\item[1.] $S_{v_0}^{(k)}$ is injective if and only if
whenever $D,D'$ are effective divisors of degree $k$ with $D \sim D'$, we
have $D = D'$.
If $S_{v_0}^{(k)}$ is injective, then $S_{v_0}^{(k')}$ is injective
for all positive integers $k' \leq k$.
\item[2.] $S_{v_0}^{(k)}$ is surjective if and only if every divisor of degree $k$ is
linearly equivalent to an effective divisor.
If $S_{v_0}^{(k)}$ is surjective, then $S_{v_0}^{(k')}$ is surjective
for all integers $k' \geq k$.
\end{itemize}
\end{lemma}

In particular, whether or not $S_{v_0}^{(k)}$ is injective (resp. surjective) is
independent of the base point $v_0$.  We therefore write $S^{(k)}$ instead of $S_{v_0}^{(k)}$ in
what follows.

 \subsection{Surjectivity of the maps $S^{(k)}$}

We recall the statement of Theorem~\ref{SurjectivityTheorem}:

\begin{theorem*}
The map $S^{(k)}$ is surjective if and only if $k \geq g$.
\end{theorem*}

\begin{proof}[Proof of Theorem~\ref{SurjectivityTheorem}]
This is an easy consequence of the Riemann-Roch theorem for graphs.
If $D$ is a divisor of degree $d \geq g$,
then since $r(K-D) \geq -1$, Riemann-Roch implies that $r(D) \geq 0$, so that $D$ is linearly equivalent to
an effective divisor.  Thus $S^{(d)}$ is surjective.
(Alternatively, we can apply (RR1) directly: if $\deg(D) \geq g$, then for all $\nu \in \N$ we have $\deg(\nu - D) < 0$ and thus
$r(\nu - D) = -1$.  By (RR1) we thus have $r(D) \geq 0$.)

Conversely, (RR1) implies that $\N \neq \emptyset$, and therefore $S^{(g-1)}$ is not surjective.
\end{proof}

\begin{remark}
This result was posed as an unsolved problem on p.~179 of \cite{BDN}.
\end{remark}

 \subsection{The chip-firing game revisited} \label{ChipFiringSection}

As mentioned earlier, Theorems \ref{ChipFiringTheorem} and \ref{SurjectivityTheorem} are equivalent.
To see this, we note the following easy lemma:

\begin{lemma}
\label{ChipFiringEquivalenceLemma}
Two divisors $D$ and $D'$ on $G$ are linearly equivalent if and only if there is a sequence of moves
in the chip firing game which transforms the configuration corresponding to $D$ into the configuration
corresponding to $D'$.
\end{lemma}

\begin{proof}
A sequence of moves in the chip-firing game can be encoded as the
function $f \in \M(G)$ for which $f(v)$ is the number of times vertex $v$ ``borrows'' a dollar
minus the number of time it ``lends'' a dollar.  (Note that the
game is ``commutative'', in the sense that the order of the moves does not matter.)
The ending configuration, starting from the initial configuration $D$
and playing the moves corresponding to $f$, is given by the divisor $D + \Delta(f)$.
So the dollar distributions achievable from the initial configuration $D$ are precisely
the divisors linearly equivalent to $D$.
\end{proof}

The equivalence between Theorem~\ref{ChipFiringTheorem} and Theorem~\ref{SurjectivityTheorem}
is now an immediate consequence of Lemma~\ref{ElementaryObservationsLemma}(1), since as we have
already noted, there is a winning strategy in the chip-firing game whose initial configuration corresponds to $D$ if and only
if $D$ is linearly equivalent to an effective divisor.
In particular, we have now proved Theorem~\ref{ChipFiringTheorem}.

 \subsection{Injectivity of the maps $S^{(k)}$}

We recall the statement of Theorem~\ref{InjectivityTheorem}.

\begin{theorem*}
The map $S^{(k)}$ is injective if and only if $G$ is
$(k+1)$-edge-connected.
\end{theorem*}

\begin{proof}
Suppose $G$ is $(k+1)$-edge-connected. Choose $v_0 \in V(G)$
arbitrarily, and let $D \in \Div^k_+(G)$. For every non-empty $A
\subseteq V(G) - \{v_0\}$, we have $\sum_{v \in A}D(v) \leq k <
\sum_{v \in A} \outdeg_A(v)$, as $\sum_{v \in A} \outdeg_A(v)$ is
equal to the size of the edge cut between $A$ and $V(G) - A$.
Therefore $D(v) < \outdeg_A(v)$ for some $v \in A$. It follows that
$D$ is $v_0$-reduced, so from Proposition~\ref{GParkingProp} we deduce
that no two distinct divisors in $\Div^k_+(G)$ are equivalent, and
therefore that the map $S^{(k)}$ is injective.

Conversely, suppose $G$ is not $(k+1)$-edge-connected. Let $C
\subseteq E(G)$ be an edge cut of size $j \leq k$, and let $X
\subseteq V(G)$ be one of the components of $G - C$. Let $D =
\sum_{v \in X} |E_v \cap C|(v)$ and  $D' = D - \Delta(\chi_X)$.
Then for each $v \in V(G)$, we have
\[
\begin{aligned}
D'(v) &= |E_v \cap C| \cdot \chi_X(v) - \sum_{e = vw \in E_v} \left( \chi_X(v) - \chi_X(w) \right) \\
&= \left\{
\begin{array}{ll}
0 & v \in X \\
|\{ e = vw \in E_v \; : \; w \in X \}| & v \not\in X \ .
\end{array}
\right.
\end{aligned}
\]
Thus $D,D' \geq 0$, $D \sim D'$, and $D \neq D'$.  It follows that
the map  $S^{(j)}$ is not injective, and consequently neither is
$S^{(k)}$.
\end{proof}

In particular, $S$ is injective if and only if every edge of $G$ is contained in a cycle.

\begin{remark}
In part (iv) of Proposition 7 in \cite{BDN}, the authors state that $S$ is injective if $G$ 
has vertex connectivity at least 2, and is not the graph consisting of one edge connecting two vertices.
However, their proof contains an error (the map $h: V \to \ZZ / n\ZZ$ which they define need not be harmonic).
In any case, Theorem~\ref{InjectivityTheorem} for $k=1$ is a more precise result.
\end{remark}

\subsection{Injectivity of the Abel-Jacobi map via circuit theory}

There is an alternate way to see that $S$ is injective if and only if $G$ is 2-edge-connected
using the theory of electrical networks (which we refer to henceforth as {\em circuit theory}).
We sketch the argument here;
see \S{15} of \cite{BiggsAPTG} for some background on electrical networks.

Consider $G$ as an electric circuit where the edges are resistors of resistance 1, and let $i_{v_0}^v(e)$ be the current flowing through
the oriented edge $e$ when one unit of current enters the circuit at $v$ and exits at $v_0$.
Let $d: C^0(G,\RR) \to C^1(G,\RR)$ and $d^*: C^1(G,\RR) \to C^0(G,\RR)$ be the usual operators on cochains
(see \S1 of \cite{BDN}).
By Kirchhoff's laws, $i_{v_0}^v$ is the unique element $i$ of $C^1(G,\RR) \cap \Image(d)$ for which $d^*(i) = (v) - (v_0)$.
It follows from the fact that $d(C^0(G,\ZZ)) = C^1(G,\ZZ)$ that
$i_{v_0}^v \in \CC^1(G,\ZZ)$ if and only if $(v) - (v_0) \in d^*(C^1(G,\ZZ)) = (d^*d)(C^0(G,\ZZ))$,
which happens if and only if $S_{v_0}(v) = 0$.

Circuit theory implies that $0 < |i_{v_0}^v(e)| \leq 1$ for every edge $e$ which belongs to a path
connecting $v$ and $v_0$.  In other words, the magnitude of the current flow is at most 1 everywhere in the circuit,
and a nonzero amount of current must flow along every path from $v$ to $v_0$.

Recall that a graph $G$ is 2-edge-connected if and only if every edge of $G$ is contained in a cycle.
So if $G$ is 2-edge-connected, then circuit theory implies that
$|i_{v_0}^v(e)| < 1$ for every edge $e$ belonging to a path connecting $v$ and $v_0$.
(Some current flows along each path from $v$ to $v_0$, and there are at least two such edge-disjoint paths.)
Therefore $i_{v_0}^v \not\in C^1(G,\ZZ)$, so $S_{v_0}(v) \neq 0$.  Since
$S_{v_0}(v) - S_{v_0}(v') = S_{v'}(v)$, this implies that $S_{v_0}$ is injective.

Conversely, if an edge $e'$ of $G$ is not contained in any cycle, then letting $v,v'$ denote the endpoints of $e'$,
it follows from circuit theory that
\[
i_{v_0}^v(e) = \left\{
\begin{array}{ll}
1 & \textrm{ if } e = e' \\
0 & \textrm{ otherwise.} \\
\end{array}
\right.
\]

Therefore $S_{v_0}(v) = S_{v_0}(v')$ and $S_{v_0}$ is not injective.

\begin{remark}
A similar argument is given in
\S{9} of \cite{Edixhoven}, although the connection with the Jacobian of a finite graph is not explicitly mentioned.
Yet another proof of the statement ``$S$ is injective if and only if $G$ is 2-edge-connected'' can be found in
Corollary 2.3 of \cite{Lorenzini1} (where the result is attributed to Hans Gerd Evertz).
\end{remark}

The circuit theory argument actually tells us something more precise about the failure of $S$ to be injective on a general
graph $G$.  Let $\overline{G}$ be the graph obtained by contracting every edge of $G$ which
is not part of a cycle, and let $\rho : G \to \overline{G}$ be the natural map.


\begin{lemma}
\label{RetractionLemma}
$\rho(v_1) = \rho(v_2)$ if and only if $(v_1) \sim (v_2)$.
\end{lemma}

\begin{proof}
$\rho(v_1) = \rho(v_2)$ if and only if there is a path from $v_1$ to $v_2$ in $G$,
none of whose edges belong to a cycle.
By circuit theory, this occurs if and only if there is a unit current flow from $v_1$ to $v_2$ which
is integral along each edge.  By the above discussion, this happens if and only if $(v_1) \sim (v_2)$.
\end{proof}

As a consequence of Lemma~\ref{RetractionLemma} and Theorem~\ref{InjectivityTheorem}, we obtain:

\begin{corollary}
\label{RetractionCor}
For every graph $G$ and every base point $v_0 \in G$, there is a commutative diagram
\[
\begin{CD}
G     @>{\rho}>>     \overline{G}     \\
@V{S}VV    @VV{\overline{S}}V     \\
\Jac(G)    @>{\rho_*}>{\cong}>   \Jac(\overline{G})   \\
\end{CD}
\]
in which $\rho_*$ is an isomorphism, $\rho$ is surjective, and
$\overline{S} = \overline{S}_{\rho(v_0)}$ is injective.
\end{corollary}


\begin{remark}

(i) It is not hard to give a rigorous proof of
Corollary~\ref{RetractionCor} which does not rely
on circuit theory by showing that the natural
map $\rho_* : \Div(G) \to \Div(G')$ given by $\rho_*(\sum a_v (v)) = \sum a_v (\rho(v))$ sends principal divisors
to principal divisors and induces a bijection $\Jac(G) \to \Jac(G')$.
We leave this as an exercise for the interested reader.

(ii) Theorem~\ref{InjectivityTheorem} and Corollary~\ref{RetractionCor} suggest that from the point of view of Abel-Jacobi theory,
the ``correct'' analogue of a Riemann surface is a 2-edge-connected graph.
This point of view resonates with the classification of Riemann surfaces by genus.
For example, there is a unique Riemann surface of genus 0 (the Riemann sphere), and there is a unique
$2$-edge-connected graph of genus 0 (the graph with one vertex and no edges).
Similarly, Riemann surfaces of genus 1 are classified up to isomorphism by a single complex number known as the ``$j$-invariant'', and
a $2$-edge-connected graph of genus 1 is isomorphic to a cycle of length $n\geq 2$, so is determined up
to isomorphism by the integer $n$.
\end{remark}

\section{Complements}
\label{ComplementsSection}

 \subsection{Morphisms between graphs}

In algebraic geometry, one is usually interested not just in Riemann surfaces themselves but
also in the holomorphic maps between them.  The most general graph-theoretic analogue of a holomorphic
map between Riemann surfaces in the context of the present paper appears to be the notion of a {\em harmonic morphism}, as defined in
\cite{Urakawa1}.  For a non-constant harmonic morphism $f : X_1 \to X_2$,
there is a graph-theoretic analogue of the classical {\em Riemann-Hurwitz formula} relating the canonical
divisor on $X_1$ to the pullback of the canonical divisor on $X_2$.
Moreover, a non-constant harmonic morphism $f : X_1 \to X_2$
induces maps $f_* : \Jac(X_1) \to \Jac(X_2)$ and $f^* : \Jac(X_2) \to \Jac(X_1)$ between the Jacobians of $X_1$ and
$X_2$ in a functorial way.
We will discuss these and other matters,
including several characterizations of ``hyperelliptic'' graphs, in a subsequent paper.

 \subsection{Generalizations}

There are some obvious ways in which one might attempt to generalize the results of this paper.  For example:

1. We have dealt in this paper only with finite unweighted graphs, but it would be interesting to generalize
our results to certain infinite graphs, as well as to weighted and/or metric graphs.



2. Can the quantity $r(D) - r(K-D)$ appearing in Theorem~\ref{GraphRiemannRochTheorem} be interpreted in a natural way
as an Euler characteristic?  In other words, is there a Serre duality theorem for graphs?

3. One could try to generalize some of the results in this paper to higher-dimensional simplicial complexes.
For example, is there a higher-dimensional generalization of Theorem~\ref{GraphRiemannRochTheorem}
analogous to the Hirzebruch-Riemann-Roch theorem in algebraic geometry?


 \subsection{Other Riemann-Roch theorems}

$\left.\right.$ \medskip

1. Metric graphs are closely related to ``tropical curves'', and in this context
Mikhalkin and Zharkov have recently announced a tropical Abel-Jacobi theorem and a tropical Riemann-Roch inequality
(see \S5.2 of \cite{MikhalkinICM}).  It appears, however, that their definition of $r(D)$ is different
from ours (this is related to the discussion in Remark~\ref{DimensionRemark}).

2. There is a Riemann-Roch formula in toric geometry having to do with lattice points and volumes of polytopes (see, e.g.,
\S5.3 of \cite{FultonToricBook}).
Our Theorem~\ref{GraphRiemannRochTheorem} appears to be of a rather different nature.

 \subsection{Connections with number theory}

The first author's original motivation for looking at the questions in this paper came from connections with
number theory.  We briefly discuss a few of these connections.

1. The Jacobian of a finite graph arises naturally in the branch of number theory known as
arithmetic geometry.  One example is the theorem of Raynaud
\cite{RaynaudPicard} relating a proper regular semistable model $\X$ for a curve $X$ over a discrete valuation ring to
the group of connected components $\Phi$ of the special fiber of the 
N{\'e}ron model of the Jacobian of $X$.
Although not usually stated in this way,
Raynaud's result essentially says that $\Phi$ is canonically isomorphic to the Jacobian of the dual graph of the special fiber of $\X$.
See \cite{Edixhoven,Lorenzini3,Lorenzini2,Lorenzini1} for further details and discussion.
Raynaud's theorem plays an important supporting role in a number of seminal papers in number theory
(see, for example, \cite{MazurMCEI} and \cite{RibetMR}).


2. The canonical divisor $K$ on a graph, as defined in (\ref{eq:canondiv}), plays a prominent role
in Zhang's refinement of Arakelov's intersection pairing on an arithmetic surface (see \cite{ZhangAP}).

3. By its definition as a ``Picard group'', the Jacobian of a finite graph $G$ can be thought of as analogous to
the ideal class group of a number field.
In particular, the number $\kappa(G)$ of spanning trees in a graph $G$,
which is the order of $\Jac(G)$, is analogous to the class number of a number field.
This analogy appears explicitly in a graph-theoretic analogue (involving the Ihara zeta function of $G$) of the analytic class number formula for the Dedekind zeta function of a number field, 
see \cite[p.11]{HortonStarkTerras}. 
See also \cite{KotaniSunada2,StarkTerrasI, StarkTerrasII,StarkTerrasIII} for further information about the
Ihara zeta function of a graph.




 \subsection{The chip-firing game of Bj{\"o}rner-Lov{\'a}sz-Shor}
 \label{BLSSection}

In this section, we describe some connections between our chip-firing game, as described in
\S\ref{ChipFiringSubsection}, and the game previously studied by Bj{\"o}rner, Lov{\'a}sz, and
Shor in \cite{BLS}.  In order to distinguish between the two, we refer to our game as the
``unconstrained chip-firing game'', and to the game from \cite{BLS} as the ``constrained chip-firing game''.

\medskip

The constrained chip-firing game is played as follows.  Each
vertex of a given (connected) graph $G$ begins with some nonnegative amount of chips, and a move consists of choosing
a vertex with at least as many chips as its degree, and having it send one chip to each of its neighbors (in which case
we say that the vertex ``fires'').
The game terminates when no vertex is able to fire.  The main results of \cite{BLS} are the following two
theorems:

\begin{theorem}[Theorem 2.1 of \cite{BLS}]
\label{BLSTheorem1}
The finiteness or non-finiteness of the constrained chip-firing game, as well as the terminal configuration and the total
number of moves when the game is finite, are independent of the particular moves made.
\end{theorem}

\begin{theorem}[Theorem 3.3 of \cite{BLS}]
\label{BLSTheorem2}
Let $N$ be the number of chips present at any point during the constrained chip-firing game.

(a) If $N > 2|E(G)| - |V(G)|$, the game is infinite.

(b) If $|E(G)| \leq N \leq 2|E(G)| - |V(G)|$, then there exists an initial configuration guaranteeing
finite termination, and also one guaranteeing an infinite game.

(c) If $N < |E(G)|$, the game terminates in a finite number of moves.
\end{theorem}

We do not have much new to say about Theorem~\ref{BLSTheorem1}.  However,
we will show that Theorem~\ref{BLSTheorem2} can be deduced from Theorem~\ref{ChipFiringTheorem}, and
conversely that Theorem~\ref{BLSTheorem2} implies the special case of Theorem~\ref{ChipFiringTheorem}
in which the initial configuration $D$ satisfies $D(v) \leq \deg(v) - 1$ for all $v \in V(G)$.

\medskip

The result which is needed to relate the two games is the following:

\begin{lemma}
\label{BorrowingLemma}
A winning strategy exists in the unconstrained chip-firing game with initial configuration $D$
if and only if there is a sequence of borrowings by vertices having a negative number of
dollars which transforms $D$ into an effective divisor.
\end{lemma}

\begin{proof}
As one direction is obvious, it suffices to show that if $D \sim E$
with $E \geq 0$, then we can get from $D$ to an effective divisor $E'$ via a (possibly empty) sequence of
borrowings by vertices having a negative number of dollars.
Since $D \sim E$, we have $E = D + \Delta(f)$ for some $f \in \M(G)$.


Let $E' = D + \Delta(f')$ be chosen so that:
\begin{itemize}
\item[(i)] $E'$ can be reached from $D$ via a (possibly empty) sequence of borrowings by vertices having a negative number of
dollars;
\item[(ii)] $f' \leq f$; and
\item[(iii)] $\sum_{v \in V(G)} f'(v)$ is maximal subject to conditions (i) and (ii).
\end{itemize}
We must have $E'(v) \geq 0$ for every $v \in V(G)$
such that $f'(v) < f(v)$, as otherwise the configuration $E' + \Delta(\chi_{ \{ v \} })$
obtained from $E'$ by having $v$ borrow a dollar from each of its neighbors
would contradict the choice of $E'$. Moreover, $E'(v) \geq E(v) \geq 0$ for every  $v \in V(G)$ such
that  $f'(v) = f(v)$. Therefore $E'$ is effective, and the lemma holds.
\end{proof}

As a consequence of Lemma~\ref{BorrowingLemma}, we can show that the two chip firing games are related by a simple
correspondence.
For $D \in \Div(G)$, define $D^\star = K^+ - D$, where
$$K^+ =\sum_{v \in V(G)}(\deg(v)-1)(v).$$
Explicitly, if $D = \sum a_v (v) \in \Div(G)$, then $D^\star = \sum a^\star_v (v)$, where
$a^\star_v = \deg(v) - 1 - a_v$.
Note that $a^\star_v \geq 0$ if and only if $a_v \leq \deg(v) - 1$, and
that $(D^\star)^\star = D$.

\begin{corollary}
\label{ChipFiringComparisonCor}
If $D = \sum a_v (v) \in \Div(G)$ with $a_v \leq \deg(v) - 1$ for all $v \in V(G)$, then
$|D| \neq \emptyset$ if and only if there is a legal sequence of firings
in the constrained chip-firing game which starts with the configuration $D^\star$
and terminates in a finite number of moves.
\end{corollary}

\begin{proof}
By Lemmas~\ref{ChipFiringEquivalenceLemma} and \ref{BorrowingLemma},
we have $|D| \neq \emptyset$ if and only if there is a sequence of
borrowings by (not necessarily distinct) vertices $v_1,\ldots,v_k$
of $G$ that leads to a nonnegative divisor $E = \sum e_v (v)$, and
such that only vertices which are in debt ever borrow. Using the
definitions, this happens if and only if firing $v_1,\ldots,v_k$ in
the constrained chip-firing game beginning at $D^\star$ yields a
legal sequence of moves ending with a divisor $E^\star = \sum
e_v^\star (v)$ having $e_v^\star \leq \deg(v) - 1$ for all $v \in
V(G)$.
\end{proof}

With the help of Corollary~\ref{ChipFiringComparisonCor},
we can use Theorem~\ref{ChipFiringTheorem} to give an
alternative proof of Theorem~\ref{BLSTheorem2}.  Indeed, suppose the constrained chip-firing game begins with
a configuration $D^\star$ with $\deg(D^\star) = N$.  Then $\deg(D) = 2|E(G)| - |V(G)| - N$, and by
Theorem~\ref{ChipFiringTheorem}, Corollary~\ref{ChipFiringComparisonCor}, and the fact that $|D'| = \emptyset$ whenever
$\deg(D') < 0$, we see that:

(a) If $\deg(D)  < 0$, the game is infinite.

(b) If $0 \leq \deg(D) \leq |E(G)| - |V(G)| = g - 1$, then there exists an initial configuration guaranteeing
finite termination, and also one guaranteeing an infinite game.

(c) If $\deg(D) > |E(G)| - |V(G)| = g - 1$, the game terminates in a finite number of moves.

This clearly implies Theorem~\ref{BLSTheorem2}.  The same reasoning shows that Theorem~\ref{BLSTheorem2} implies
Theorem~\ref{ChipFiringTheorem} in the special case where $D(v) \leq \deg(v) - 1$ for all $v \in V(G)$.

\begin{remark}

(i) Theorem 1 of \cite{Thorup} provides a short and elegant proof of Theorem~\ref{BLSTheorem1}, and can also be used to show
that in the unconstrained chip-firing game with initial configuration $D$, every sequence of borrowings from vertices having a
negative number of dollars is either infinite (if $|D| = \emptyset$)
or else terminates in the same number of moves (when $|D| \neq \emptyset$).  In the latter case,
just as in the constrained chip-firing game, the terminal configuration is independent of the particular moves made.

(ii) If any (or equivalently, every) sequence of borrowings by
vertices in debt starting with the initial configuration $D$
terminates, then by an argument from \cite{Tardos} it terminates in
at most $\deg^+(D)d(G)|V(G)|$ steps, where $d(G)$ denotes the
diameter of $G$, i.e., the maximum path-distance between two
vertices of $G$. Thus there exists an algorithm for determining
whether $|D| = \emptyset$ whose running time is bounded from above
by $\deg^+(D)d(G)|V(G)|$.

\end{remark}

 \subsection{Reduced divisors and critical configurations}
 \label{BiggsSection}

In ~\cite{BiggsCF,BiggsTutte} (see also Chapter 14 of \cite{GR}), Biggs studies the {\em critical group} of a graph, which he
defines in terms of a certain chip-firing game played on the vertices of the graph.
One of Biggs' results is that the critical group is isomorphic to $\Jac(G)$.
In this section, we describe a one-to-one correspondence between elements of Biggs' critical
group and $v_0$-reduced divisors, as defined in \S\ref{GParkingSection}.
In order to do this, we first need to translate Biggs' definitions into the language of divisors.


Let $v_0 \in V(G)$ be a fixed base vertex, and let $v_1,v_2, \ldots,
v_{n-1}$ be an ordering of the vertices in $V(G) - \{v_0\}$, where
$n=|V(G)|$. We say that a divisor $D$ is \emph{$v_0$-critical with
respect to the ordering $v_1,v_2, \ldots, v_{n-1}$} if  for every $v
\in V(G) - \{v_0\}$ we have $0 \le D(v) < \deg(v)-1$, and if for
every $1 \leq k \leq n-1$ we have $D_k(v) \geq 0$, where $$D_k=D -
\sum_{i=0}^k \Delta(\chi_{\{v_i\}}).$$ We say that a divisor $D$ is
\emph{$v_0$-critical} if it is $v_0$-critical with respect to some
ordering of $V(G) - \{v_0\}$.

We remark on some technical differences between the above definition
and the definition given in ~\cite{BiggsCF}. In ~\cite{BiggsCF}, only configurations
for which the total amount of money is zero are considered.
Also, the definition of a critical configuration given in ~\cite{BiggsCF},
when translated directly into the language of divisors,
would appear to be slightly different from ours; however,
the two definitions are in fact equivalent by Lemma 2.6 of
~\cite{BiggsCF}.

It follows from the results of ~\cite{BiggsCF} and \cite{BiggsTutte} that given $v_0 \in
V(G)$, every equivalence class of $\Div(G)$ contains a unique
$v_0$-critical divisor. This observation suggests a
relationship between $v_0$-reduced and $v_0$-critical divisors.
In the following lemma, we show that in fact there exists a natural bijection
between the two.

\begin{lemma}
\label{ReducedCriticalLemma}
A divisor $D$ is $v_0$-reduced if and only if the
divisor $D^\star = K^+ - D$ is $v_0$-critical.
\end{lemma}

\begin{proof}
Let $n=|V(G)|$. Suppose that $D$ is $v_0$-reduced, and
define $v_1,v_2, \ldots,v_{n-1}$ as in the proof of Theorem \ref{MainGraphTheorem}.
We claim that $D^\star$ is $v_0$-critical with respect to this ordering
of $V(G)  - \{v_0\}$.

Write $D^\star_k$ for $(D^\star)_k = K^+ - D - \Delta(\chi_{B_k})$, where
$B_k=\{v_0,v_1,\ldots,v_{k}\}$. Let $v \in V(G) - \{v_0\}$. We have
$0 \le D(v) <\outdeg_{\{v\}}(v)= \deg(v)$ and therefore $0 \le
D^\star(v) < \deg(v)$. It remains to prove that $0 \le D^\star_k(v)$
for every $1 \leq k \leq n-1$.

If $v \not \in B_k$, then $D^\star_k(v) \geq D^\star(v) \geq 0$. Otherwise
$v=v_l$ for some $0< l \leq k$, and
\begin{align*}
D^\star_k(v_l) \geq D^\star_{l}(v_l) &= \deg(v_l)- 1 - D(v_l) -
\outdeg_{B_l}(v_l) \\
&= \left( \deg(v_l) - \outdeg_{B_l}(v_l) \right) - D(v_l) -1 \\
&= \outdeg_{V(G)-B_{l-1}}(v_l)-D(v_l)-1 \\
&\geq 0 \ ,
\end{align*}
where the last inequality follows from the definition of $v_l$.
(Here we have used the fact that if $A \subseteq V(G)$, then
$\outdeg_A(v) + \outdeg_{V(G) - A}(v) = \deg(v)$ for all $v \in V(G)$, and if $v \in A$, then
$\outdeg_A(v) = \outdeg_{A - \{ v \} }(v)$.)
It follows that $D^\star$ is $v_0$-critical with respect to the given order,
as desired.

\medskip

Now suppose $D^\star$ is $v_0$-critical with respect to the ordering
$v_1$,$v_2, \ldots$, $v_{n-1}$. Consider a non-empty subset $A \subseteq V(G)
- \{v_0\}$, and let $v_l$ be the vertex in $A$ having the smallest index. We
have
\[
0 \leq D^\star_l(v_l)= \outdeg_{V(G)-B_{l-1}}(v_l)-D(v_l)-1 \ ,
\]
where $B_{l-1}$ is defined as above. Moreover, $B_{l-1} \cap A =
\emptyset$, and therefore
\[
D(v_l) < \outdeg_{V(G)-B_{l-1}}(v_l) \leq \outdeg_{A}(v_l) \ .
\]
As $A \subseteq V(G) - \{v_0\}$ was arbitrary, we
conclude that $D$ is $v_0$-reduced.
\end{proof}

\begin{remark}
Lemma~\ref{ReducedCriticalLemma} explains some of the parallels found in
the literature between certain results concerning $G$-parking functions and critical configurations.
As two examples, we mention:
\begin{itemize}
\item[(i)] The construction of explicit bijections between $G$-parking functions and spanning trees
from \cite{CP}, and between critical configurations and spanning trees in
\cite{BiggsWinkler}.
\item[(ii)] The relationship between $G$-parking functions and
the Tutte polynomial, as described in \cite{PC}, and between critical configurations and
the Tutte polynomial, as described in \cite{Merino2} and \cite{BiggsTutte}.
\end{itemize}
\end{remark}



\appendix
\section{Riemann surfaces and their Jacobians}
\label{RSAppendix}

The theory of Riemann surfaces and their Jacobians is
one of the major accomplishments of $19^{\rm th}$ century
mathematics, and it continues to this day to have significant applications.
We cannot hope to give the reader a complete overview
of this vast subject, so we will just touch on a few of the highlights of the theory in order to draw out the connections
with graph theory.  We recommend \cite{Miranda} as a good introduction to the theory of Riemann surfaces and their Jacobians;
see also \cite{ACGH,FarkasKra,GH,Hartshorne,Mumford,Murty}.

\medskip

A (compact) {\em Riemann surface} $X$ is a one-dimensional connected complex manifold, i.e., a two-dimensional connected compact real
manifold endowed with a maximal atlas $\{ U_\alpha, z_\alpha \}$ for which the transition functions
\[
f_{\alpha \beta} = z_{\alpha} \circ z_{\beta}^{-1} : z_{\beta}(U_\alpha \cap U_\beta) \to
z_{\alpha}(U_\alpha \cap U_\beta)
\]
are holomorphic whenever $U_\alpha \cap U_\beta \neq \emptyset$.

The simplest example of a Riemann surface is the {\em Riemann sphere} $\CC \cup \{\infty \}$.

Since a Riemann surface looks locally like an open subset of $\CC$,
there is a natural notion of what is means for a function $f : X \to \CC$ (resp.
$f : X \to \CC \cup \{ \infty \}$) to be {\em holomorphic} (resp. {\rm meromorphic}):
we say that $f$ is holomorphic (resp. meromorphic) if $f \circ z^{-1}$ is holomorphic (resp. meromorphic)
for every coordinate chart $(U,z)$.

\medskip

A {\em $1$-form} $\omega$ on a Riemann surface $X$ is a collection of $1$-forms $\omega_x dx + \omega_y dy$
on each coordinate chart $(U,z)$ (where $z = x+iy$) satisfying suitable compatibility relations on
overlapping charts.  A $1$-form is {\em holomorphic} if $\omega_x$ and $\omega_y$ are holomorphic and
$\omega_y = i\omega_x$.  Locally, every holomorphic $1$-form is equal to $f(z)dz$ with $f$ a holomorphic
function.  Finally, a $1$-form is {\em meromorphic} if it is holomorphic outside a finite set of points
and can be represented locally as $f(z)dz$ with $f$ a meromorphic function.

\medskip

Riemann surfaces are classified by a nonnegative integer $g$ called the {\em genus}.  There are several equivalent
characterizations of the genus of a Riemann surface; for example, $2g$ is the topological genus of $X$,
i.e., $\dim_\RR H_1(X,\RR)$, and $g$ is the complex dimension of the space of
holomorphic $1$-forms on $X$.  A Riemann surface has genus 0 if and only if it is isomorphic to the Riemann sphere.





\medskip

Let $\Div(X)$ be the free abelian group on the set of vertices of $X$;
elements of $\Div(X)$ are called {\em divisors} on $X$ and are usually written as
$\sum_{P \in X} a_P (P)$, where each $a_P$ is an integer and all but finitely many of the $a_P$'s are zero.
A divisor $E \in \Div(X)$ is called {\em effective} if
$E \geq 0$.

There is a natural {\em degree function} $\deg : \Div(X) \to \ZZ$ given for $D = \sum a_P (P)$ by
\[
\deg(D) = \sum_{P \in X} a_P \ .
\]

If $\M(X)$ denotes the space of meromorphic functions on $X$, then for every nonzero $f \in \M(X)$ and every $P \in X$,
one can define, using local coordinates, the {\em order of vanishing} $\ord_P(f)$ of $f$ at $P$.
For all but finitely many $P \in X$, one has $\ord_P(f) = 0$.
The {\em divisor of $f$} is then defined to be
\begin{equation}
\label{RSDivDef}
{\rm div}(f) = \sum_{P \in X} \ord_P(f) (P) \ .
\end{equation}
The divisor of a nonzero meromorphic function $f$ is called a {\em principal divisor}.
A fundamental fact about Riemann surfaces is that $\deg({\rm div}(f)) = 0$, which means that
$f$ has the same number of zeros as poles (counting multiplicities).
Therefore $\Prin(X)$ (the set of all principal divisors) is a subgroup of the group
$\Div^0(X)$ of divisors of degree zero.

\medskip

The Jacobian $\Jac(X)$ of $X$ (also denoted $\Pic^0(X)$) is defined to be the quotient group
\begin{equation}
\label{RSJacobianDef}
\Jac(X) = \frac{\Div^0(X)}{\Prin(X)} \ .
\end{equation}
The abelian group $\Jac(X)$ is naturally endowed with the structure of a (projective)
compact complex manifold of dimension $g$,
i.e., $\Jac(X)$ is an {\em abelian variety}.

Two divisors $D,D'$ on $X$ are called {\em linearly equivalent} if their difference is a principal divisor.
Thus $\Jac(X)$ classifies the degree zero divisors on $X$ modulo linear equivalence.





\medskip

If we fix a base point $P_0 \in \Jac(X)$, we can define the {\em Abel-Jacobi map} $S_{P_0} : X \to \Jac(X)$
by the formula
\begin{equation}
\label{eq:RSAbelJacobiMap}
S_{P_0}(P) = [(P) - (P_0)] \ ,
\end{equation}
where $[D]$ denotes the class in $\Jac(X)$ of $D \in \Div^0(X)$.
We write $S$ instead of $S_{P_0}$ when the base point $P_0$ is understood.

We can also define, for each $k \geq 1$, the map $S_{P_0}^{(k)} : \Div_+^k(X) \to \Jac(X)$
by
\[
S_{P_0}^{(k)}((P_1) + \cdots + (P_k)) = S_{P_0}(P_1) + S_{P_0}(P_2) + \cdots + S_{P_0}(P_k) \ ,
\]
where $\Div_+^k(X)$ denotes the set of effective divisors of degree $k$ on $X$.

The map $S_{v_0}$ can be characterized by the following universal property:
If $\varphi$ is a holomorphic map from $X$ to an abelian variety $A$ taking $P_0$ to $0$,
then there is a unique homomorphism
$\psi : \Jac(X) \to A$ such that $\varphi = \psi \circ S_{P_0}$.

A classical result about the maps $S^{(k)}$ is the following:

\begin{theorem}
\label{RSSurjectivityTheorem}
$S^{(k)}$ is surjective if and only if $k \geq g$.
\end{theorem}

The surjectivity of $S^{(g)}$ is usually referred to as {\em Jacobi's inversion theorem};
it is equivalent to the statement that every divisor of degree at least $g$ on $X$ is linearly equivalent to
an effective divisor.

Another classical fact is:

\begin{theorem}
\label{RSInjectivityTheorem}
The Abel-Jacobi map $S$ is injective if and only if $g \geq 1$.
\end{theorem}

\medskip

Let $D$ be a divisor on $X$.  The {\em linear system associated to $D$} is defined to
be the set $|D|$ of all effective divisors linearly equivalent to $D$:
\[
|D| = \{ E \in \Div(X) \; : \; E \geq 0, \; E \sim D \} \ .
\]

\medskip

The {\em dimension} $r(D)$ of the linear system $|D|$ is defined to be one less than the dimension of
$L(D)$, where
\[
L(D) = \{ f \in \M(X) \; : \; \divisor(f) \geq -D \}
\]
is the finite-dimensional $\CC$-vector space consisting of all meromorphic functions for which
$\divisor(f) + D$ is effective.  There is a natural identification
\[
|D| = \left( L(D) - \{ 0 \} \right) / \CC^*
\]
of $|D|$ with the projectivization of $L(D)$.
It is easy to see that $r(D)$ depends only on the linear equivalence class of $D$.

\begin{remark}
\label{RSDimensionRemark}
In the graph-theoretic setting, the analogue of $L(D)$ is no longer a vector space.
Therefore it is useful to have a more intrinsic characterization of the quantity $r(D)$ in terms
of $|D|$ only.  Such a characterization is in fact well-known
(see, e.g., p.250 of \cite{GH} or \S{III.8.15} of \cite{FarkasKra}):
$r(D) \geq -1$ for all $D$, and for each $s \geq 0$ we have
$r(D) \geq s$ if and only if $|D - E| \neq \emptyset$
for all effective divisors $E$ of degree $s$.
\end{remark}

\medskip

Given a nonzero meromorphic $1$-form $\omega$ on $X$, one can define (using local coordinates) the order of vanishing of $\omega$ at
a point $P \in X$, and the divisor $\divisor(\omega)$ of $\omega$ is then defined as in (\ref{RSDivDef}).
The degree of $\divisor(\omega)$ is $2g-2$ for every $\omega$, and if $\omega,\omega'$ are both
nonzero meromorphic $1$-forms on $X$, the quotient $\omega / \omega'$ is a nonzero meromorphic function on $X$, and
thus $\divisor(\omega)$ and $\divisor(\omega')$ are linearly equivalent.

The {\em canonical divisor class} $K_X$ on $X$ is defined to be the linear equivalence class of $\divisor(\omega)$ for any
nonzero meromorphic $1$-form $\omega$.





\medskip

The following result, known as the {\em Riemann-Roch theorem}, is widely regarded as
the single most important result in the theory of Riemann surfaces.

\begin{theorem}[Riemann-Roch]
\label{RiemannRochTheorem}
Let $X$ be a Riemann surface with canonical divisor class $K$, and let $D$ be a divisor on $X$.  Then
\[
r(D) - r(K - D) = \deg(D) + 1 - g \ .
\]
\end{theorem}


The importance of Theorem~\ref{RiemannRochTheorem} stems from the large number of
applications which it has; see, e.g., Chapters VI and VII of \cite{Miranda} and
Chapter IV of \cite{Hartshorne}.


\medskip

Finally, we discuss {\em Abel's theorem}, which gives an alternative characterization of $\Jac(X)$
and the Abel-Jacobi map $S_{P_0} : X \to \Jac(X)$.

Choose a base point $P_0 \in X$, and
let $\Omega^1(X)$ denote the space of holomorphic $1$-forms on $X$.
Every (integral) homology class $\gamma \in H_1(X,\ZZ)$ defines an element $\int_{\gamma}$ of
the dual space $\Omega^1(X)^*$ via integration:
\[
\int_{\gamma} : \omega \mapsto \int_{\gamma} \omega \in \CC \ .
\]
A linear functional $\lambda : \Omega^1(X) \to \CC$ is called a {\em period} if it is of the form
$\int_{\gamma}$ for some $\gamma \in H_1(X,\ZZ)$.  We let $\Lambda$ denote the set of periods; it is a lattice
in $\Omega^1(X)^*$.

For each point $P \in X$, choose a path $\gamma_P$ in $X$ from $P_0$ to $P$, and define
$A_{P_0} : X \to \Omega^1(X)^* / \Lambda$ by sending $P$
to class of the linear functional $\int_{\gamma_P}$ given by integration along $\gamma_P$.
This is well-defined, since if $\gamma'_P$ is another path from $P_0$ to $P$, then
the 1-chain $\gamma_P - \gamma'_P$ is closed and therefore defines an integral homology class.

We can extend the map $A_{P_0}$ by linearity to a homomorphism from $\Div(X)$ to $\Omega^1(X)^* / \Lambda$.
Restricting to $\Div^0(X)$ gives a canonical map $A : \Div^0(X) \to \Omega^1(X)^* / \Lambda$
which does not depend on the choice of base point $P_0$.

\begin{theorem}[Abel's Theorem]
\label{RSAbelTheorem}
The map $A$ is surjective, and its kernel is precisely $\Prin(X)$.  Therefore $A$ induces an isomorphism of $\Jac(X)$
onto $\Omega^1(X)^* / \Lambda$.
Moreover, we have $A_{P_0} = A \circ S_{P_0}$, i.e., $A_{P_0}$ coincides with
the Abel-Jacobi map $S_{P_0}$ under the identification of
$\Jac(X)$ and $\Omega^1(X)^* / \Lambda$ furnished by $A$.
\end{theorem}

In particular, if $D$ is a divisor of degree zero on $X$, then $D$ is the divisor of a meromorphic function on $X$
if and only if $A(D) = 0$.


\section{Abel's theorem for graphs}
\label{GraphAbelTheoremSection}

For the sake of completeness, we recall from \cite{BDN}
a graph-theoretic analogue of Abel's theorem (Theorem~\ref{RSAbelTheorem}).
See also \cite{Nagnibeda} and \S{28-29} of \cite{BiggsAPTG} for further details.

\medskip

Choose an orientation of the graph $G$, i.e., for each edge $e$ pick a vertex $e_{+}$ incident to $e$,
and let $e_{-}$ be the other vertex incident to $e$.
Let $C^0(G,\RR)$ be the $\RR$-vector space consisting of all functions $f : V(G) \to \RR$.
Inside this space, we have the lattice $C^0(G,\ZZ)$ consisting of the integer valued functions.
Similarly, we can consider the space $C^1(G,\RR)$ of all functions $h : E(G) \to \RR$
and the corresponding lattice $C^1(G,\ZZ)$.  We equip $C^0(G,\RR)$ and $C^1(G,\RR)$
with the inner products given by
\begin{equation}
\label{eq:IP1}
\langle f_1, f_2 \rangle = \sum_{v \in V(G)} f_1(v) f_2(v)
\end{equation}
and
\begin{equation}
\label{eq:IP2}
\langle h_1, h_2 \rangle = \sum_{e \in E(G)} h_1(e) h_2(e) \ .
\end{equation}

Define the {\em exterior differential} $d : C^0(G, \RR) \to C^1(G, \RR)$ by the formula
\[
df(e) = f(e_{+}) - f(e_{-}) \ .
\]
The adjoint $d^* : C^1(G, \RR) \to C^0(G, \RR)$ of $d$ with respect to the inner products
(\ref{eq:IP1}) and ({\ref{eq:IP2}) is given by
\[
(d^*h)(v) = \sum_{\substack{e \in E(G) \\ e_{+} = v }} h(e) - \sum_{\substack{e \in E(G) \\ e_{-} = v }} h(e) \ .
\]
It is easily checked that $\Delta = d^* d : C^0(G, \RR) \to \C^0(G, \RR)$ is independent of the choice of orientation, and
can be identified with the Laplacian operator on $G$, i.e.:
\[
(d^*df)(v) = \deg(v) f(v) - \sum_{e = wv \in E_v} f(w) \ .
\]

There is an orthogonal decomposition
\[
C^1(G, \RR) = \Ker(d^*) \oplus \Image(d) \ ,
\]
where $\Ker(d^*)$ is the {\em flow space} (or {\em cycle space}) and
$\Image(d)$ is the {\em cut space} (or {\em potential space}).


\medskip

The {\em lattice of integral flows} is defined to be $\Lambda^1(G) = \Ker(d^*) \cap C^1(G,\ZZ)$,
and the {\em lattice of integral cuts} is defined to be
$N^1(G) = \Image(d) \cap C^1(G, \ZZ)$.

For a lattice $\Lambda$ in a Euclidean inner product space $V$, the {\em dual lattice} $\Lambda^{\#}$ is defined to be
\[
\Lambda^{\#} = \{ x \in V \; : \; \langle x, \lambda \rangle \in \ZZ \textrm{ for all } \lambda \in \Lambda \} \ .
\]

A lattice $\Lambda$ is called {\em integral} if $\langle \lambda, \mu \rangle \in \ZZ$ for all $\lambda, \mu \in \Lambda$; this
is equivalent to requiring that $\Lambda \subseteq \Lambda^{\#}$.
Clearly $\Lambda^1(G)$ and $N^1(G)$ are integral lattices.

In the statement of the following theorem, the dual of $\Lambda^1(G)$ (resp. $\N^1(G)$) is defined with respect to
the ambient space $\Ker(d^*)$ (resp. $\Image(d)$).

\begin{theorem}
\label{MultiIsomorphismTheorem}
The groups $C^1(G,\ZZ) / (\Lambda^1(G) \oplus N^1(G)), \Lambda^1(G)^{\#} / \Lambda^1(G)$, and $N^1(G)^{\#} / N^1(G)$
are all isomorphic.
\end{theorem}

Choose a base vertex $v_0 \in G$.  One can describe a map $A_{v_0} : G \to J(G) := \Lambda^1(G)^{\#} / \Lambda^1(G)$
as follows.  For any $v \in V(G)$, choose a path $\gamma$ from $v_0$ to $v$, which may be identified
in the obvious way with an element of $C^1(G,\ZZ)$.
If $\gamma'$ is any other path from $v_0$ to $v$, then $\gamma - \gamma' \in \Lambda^1(G)$.
Since $\langle \gamma, \lambda \rangle \in \ZZ$ for every $\lambda \in \Lambda^1(G)$,
$\gamma$ determines an element $A_{\gamma}$ of $\Lambda^1(G)^{\#}$.  We define
$A_{v_0}(v)$ to be the class of $A_{\gamma}$ in $\Lambda^1(G)^{\#} / \Lambda^1(G)$; this is independent of the choice of $\gamma$.

We can extend the map $A_{v_0}$ by linearity to a homomorphism from $\Div(G)$ to $\Lambda^1(G)^{\#} / \Lambda^1(G)$.
Restricting to $\Div^0(G)$ gives a canonical map $A : \Div^0(G) \to J(G)$ which does not depend on
the choice of base point $v_0$.

\begin{theorem}[Abel's Theorem for Graphs]
\label{GraphAbelTheorem}
The map $A$ is surjective, and its kernel is precisely $\Prin(G)$.  Therefore $A$ induces an isomorphism of $\Jac(G)$
onto $\Lambda^1(G)^{\#} / \Lambda^1(G)$.  Moreover, we have $A_{v_0} = A \circ S_{v_0}$, i.e., $A_{v_0}$ coincides with
the Abel-Jacobi map $S_{v_0}$ defined by (\ref{eq:AbelJacobiMap}) under the identification of
$\Jac(G)$ and $J(G)$ furnished by $A$.
\end{theorem}

Consequently, if $D$ is a divisor of degree zero on $G$, then $D$ is principal
if and only if $A(D) = 0$.
For proofs of Theorems~\ref{MultiIsomorphismTheorem} and \ref{GraphAbelTheorem},
see \cite{BDN} and \S{24-29} of \cite{BiggsAPTG}.

\begin{remark}
The lattices $\Lambda^1(G)$ and $N^1(G)$ have a number of interesting combinatorial properties.
For example, it is shown in Propositions 1 and 2 of \cite{BDN} that $\Lambda^1(G)$ is even if and only if $G$ is bipartite,
and $N^1(G)$ is even if and only if $G$ is Eulerian.
Moreover, the length of the shortest nonzero vector in $\Lambda^1(G)$ is the girth
of $G$, and the length of the shortest nonzero vector in $N^1(G)$ is the edge connectivity of $G$.
And of course, it follows from Theorem~\ref{MultiIsomorphismTheorem} that both
$|\Lambda^1(G)^{\#} / \Lambda^1(G)|$ and $|N^1(G)^{\#} / N^1(G)|$ are equal to the number of spanning trees in $G$.
\end{remark}



\bibliographystyle{plain}
\bibliography{graphs}

\end{document}